\documentclass{ceb-l}

\usepackage{amsthm}

\usepackage{amssymb}

\usepackage{graphicx}

\usepackage{pinlabel}

\newtheorem{theorem}{Theorem}[section]

\theoremstyle{definition}

\newtheorem{example}[theorem]{Example}

\theoremstyle{plain} 
\newcommand{\thistheoremname}{}
\newtheorem{genericthm}[theorem]{\thistheoremname}

\newtheorem*{genericthm*}{\thistheoremname}
\newenvironment{namedthm*}[1]
  {\renewcommand{\thistheoremname}{#1}%
   \begin{genericthm*}}
  {\end{genericthm*}}

\theoremstyle{remark}

\numberwithin{equation}{section}

\newcommand{\Z}{\mathbb{Z}}
\newcommand{\R}{\mathbb{R}}

\newcommand\Cone{\operatorname{Cone}}

\begin{document}

\title{Getting a handle on the Conway knot}


\author{Jennifer Hom}
\thanks{The author was partially supported by NSF grants DMS-1552285 and DMS-2104144. This article is associated with a lecture given by the author in the Current Events Bulletin session of the Joint Mathematics Meetings in January 2021.}
\address{School of Mathematics, Georgia Institute of Technology, Atlanta, GA, USA}
\email{hom@math.gatech.edu}
\thanks{}



\date{}

\dedicatory{}

\begin{abstract}
A knot is said to be slice if it bounds a smooth disk in the 4-ball. For 50 years, it was unknown whether a certain 11 crossing knot, called the Conway knot, was slice or not, and until recently, this was the only one of the thousands of knots with fewer than 13 crossings whose slice-status remained a mystery. We will describe Lisa Piccirillo’s proof that the Conway knot is not slice. The main idea of her proof is given in the title of this article.
\end{abstract}

\maketitle


\section{Introduction}

Here is a 3-ball: 

\begin{center}
\includegraphics[scale=0.65]{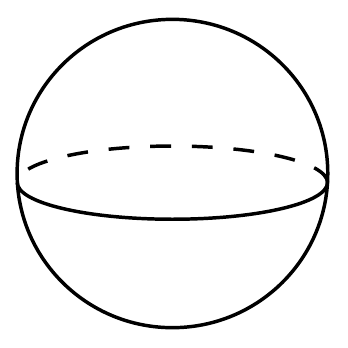} 
\end{center}

\noindent and here is a 3-ball with a handle attached:

\begin{center}
\includegraphics[scale=0.65]{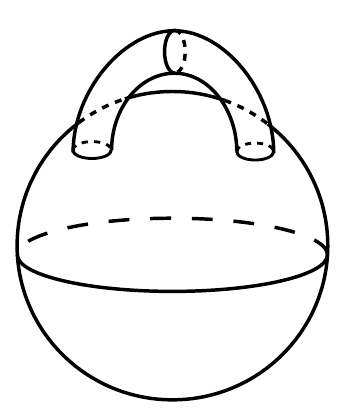} 
\end{center}

\noindent This is the Conway knot:

\begin{center}
\includegraphics[scale=0.8]{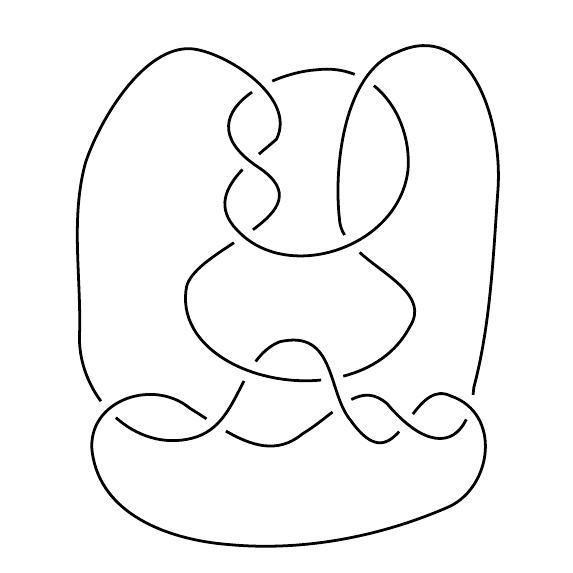} 
\end{center}

Our knots will live in the 3-sphere $S^3$, which is the boundary of the 4-ball $B^4$. A knot is \emph{slice} if it bounds a smooth disk in the 4-ball. The term slice comes from the fact that such knots are cross sections (i.e., slices) of higher dimensional knots.

\begin{namedthm*}{Main Theorem}[Piccirillo \cite{Piccirillo}]
The Conway knot is not slice.
\end{namedthm*}

This knot is not the Conway knot:
\begin{center}
\includegraphics[scale=0.8]{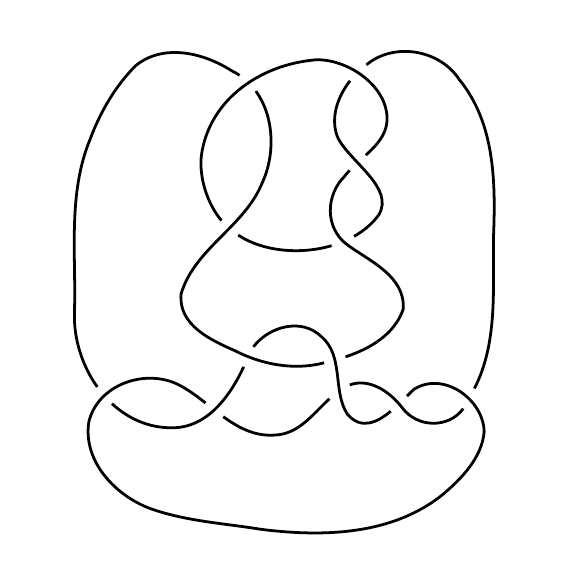} 
\end{center}
It is called the Kinoshita-Terasaka knot, and it is related to the Conway knot by \emph{mutation}, that is, we cut out a ball containing part of the knot, rotate it $180^\circ$, and glue it back in.
\begin{center}
\labellist
	\pinlabel {$180^\circ$} at 166 148
\endlabellist
\includegraphics[scale=0.8]{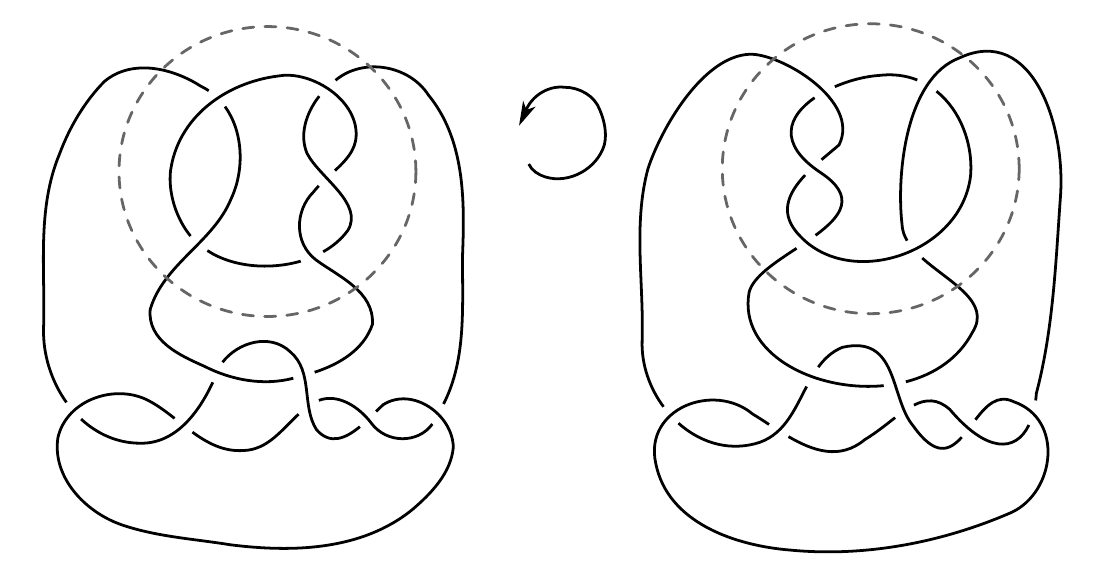} 
\end{center}

The Kinoshita-Terasaka knot is slice. Here is a slightly different diagram of the Kinoshita-Terasaka knot. As we can see, it bounds an immersed disk in $S^3$:
\begin{center}
\includegraphics[scale=0.8]{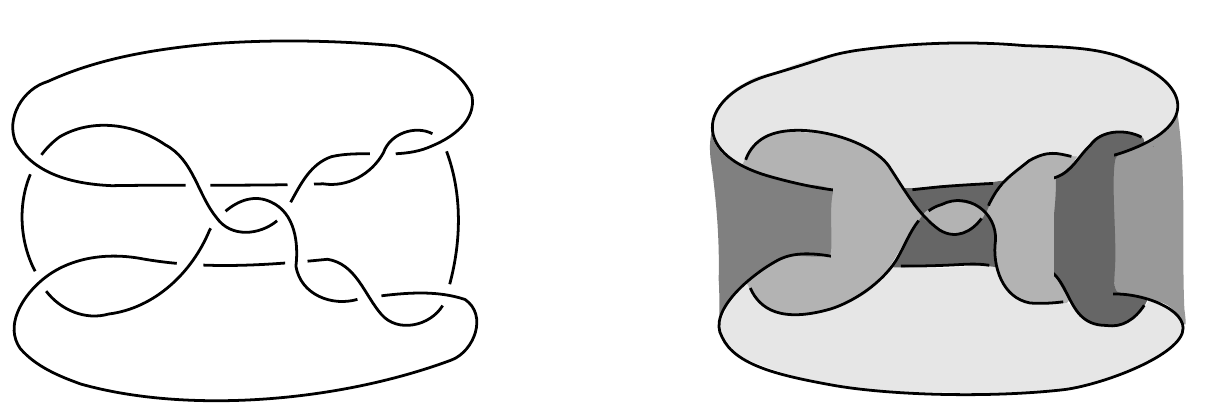} 
\end{center}
Thinking of this immersed disk as sitting in the $S^3$ boundary of the 4-ball, we can push the surface into the 4-ball and eliminate the arcs of self-intersection by pushing one sheet of the surface near the arc deeper into the 4-ball, giving us an embedded disk in the 4-ball.

One way to study knots is to use a knot invariant. A knot invariant is a mathematical object (like a number, a polynomial, or a group) that we assign to a knot. Knot invariants can be used to distinguish knots. Certain knot invariants obstruct a knot from being slice. One such invariant is Rasmussen's $s$-invariant, which to a knot $K$ assigns an integer $s(K)$. If $s(K) \neq 0$, then $K$ is not slice. 

Since the Conway knot and the Kinoshita-Terasaka knots are mutants, they have a lot in common. For example, the $s$-invariant of both knots is zero. In fact, all known knot invariants that obstruct sliceness vanish for the Conway knot. That leads one to wonder: how did Piccirillo show that the Conway knot is not slice? Her key idea was to find some other knot $K'$ such that the Conway knot is slice if and only if $K'$ is slice, and to obstruct $K'$ from being slice. The goal of this article is to give some context for her result and sketch the main ideas of her proof.

\section{Telling knots apart}\label{sec:tellingapart}

The fundamental group is one of the first algebraic invariants encountered in a topology class. A knot is homeomorphic to $S^1$, so its fundamental group is always isomorphic to the integers. However, instead of studying the knot, we can study the space around the knot. That is, we consider the \emph{knot complement}, consisting of the 3-sphere minus a neighborhood of the knot. The \emph{knot group} is the fundamental group of the knot complement.

Typically, one studies knots up to \emph{ambient isotopy}. Intuitively, this means that we can wiggle and stretch our knot, but we cannot cut it nor let it pass through itself. Since isotopic knots have homeomorphic complements and homeomorphic spaces have isomorphic fundamental groups, the knot group is an invariant of the isotopy class of a knot.

Here are two knots, the unknot and the trefoil:
\begin{center}
\includegraphics[scale=0.8]{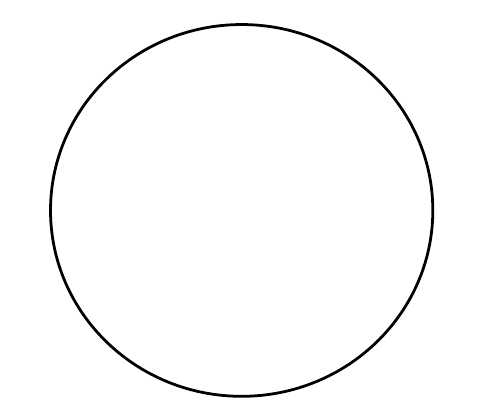} 
\includegraphics[scale=0.8]{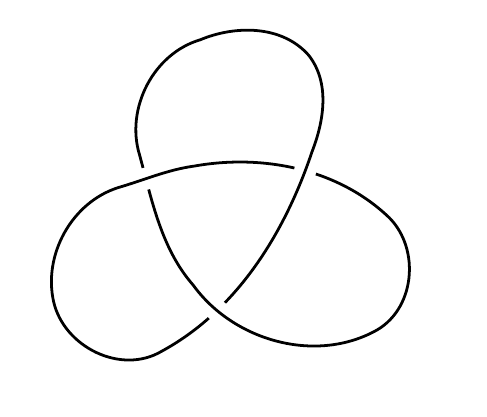} 
\end{center}

\begin{example}
The knot group of the unknot is $\Z$.

\end{example}

\begin{example}
The knot group of the trefoil is $\langle x, y \mid x^2=y^3 \rangle$. This group is non-abelian, since it surjects onto the symmetric group $S_3$. Therefore, the trefoil and the unknot are different.
\end{example}

\noindent Riley \cite{Riley} distinguished the Kinoshita-Terasaka knot and the Conway knot up to isotopy, using a delicate analysis of their fundamental groups.

Since it can often be difficult to tell if two group presentations describe isomorphic groups, it can be convenient to pass to more tractable invariants. One example is the Alexander polynomial, denoted $\Delta(t)$, which Fox \cite{FoxI} showed can be algorithmically computed from a group presentation for the knot complement.

\begin{example}
The Alexander polynomial of the unknot is $1$.
\end{example}

\begin{example}
The Alexander polynomial of the trefoil is $t^2-t+1$.
\end{example}

\begin{example}
The Conway knot and the Kinoshita-Terasaka knot both have Alexander polynomial $1$.
\end{example}

The Alexander polynomial is invariant under mutation, which explains why the Conway knot and the Kinoshita-Terasaka knot have the same Alexander polynomial. There are several other polynomial knots invariants, such as the Jones, HOMFLY-PT, and Kauffman polynomials, all of which are also invariant under mutation. Knot Floer homology \cite{OSknots} and Khovanov homology \cite{Khovanov} categorify the Alexander and Jones polynomials; that is, to a knot, they assign a graded vector space whose graded Euler characteristic is the desired polynomial. A certain version of knot Floer homology is invariant under mutation \cite{Zibrowius}, as are versions of Khovanov homology \cite{Bloom, Wehrli}. Moreover, Rasmussen's $s$-invariant is invariant under mutation \cite{KWZ}; this gives a quick way to determine that the $s$-invariant of the Conway knot is zero, since it is the mutant of a slice knot.



As we already observed, isotopic knots have homeomorphic complements. What about the converse? If two knots have 
homeomorphic complements, then are they isotopic? This question was answered in the affirmative in 1989 by Cameron Gordon and John Luecke \cite{GordonLuecke}, who proved that knots are determined by their complements. 
This is in contrast to links. 
For example, the two links below have homeomorphic complements, but are not isotopic, since in the first, both components are unknots, while in the second, one component is the trefoil.
\begin{center}
\includegraphics[scale=0.8]{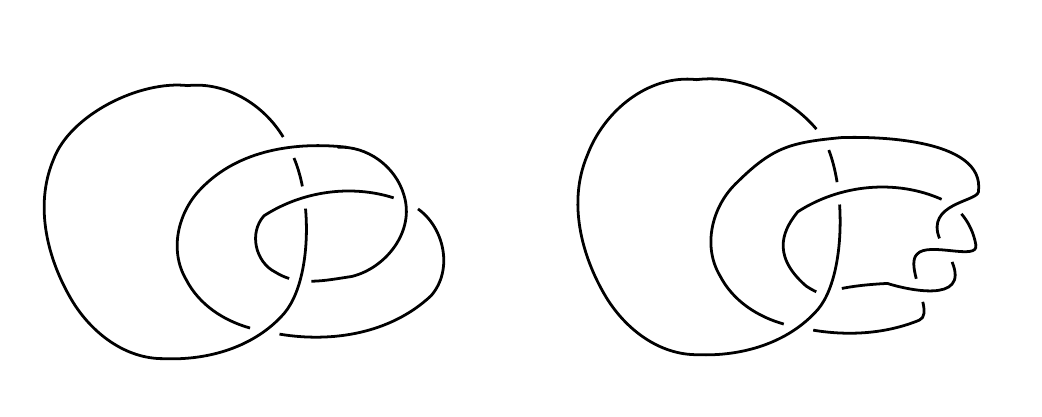} 
\end{center}

\section{Measuring the complexity of a knot}

How can we measure the complexity of a knot $K$? One such measure is the \emph{unknotting number}, denoted $u(K)$, which is the minimal number of times a knot must be passed through itself to untie it. Both the Conway knot and the Kinoshita-Terasaka knot can be unknotted by changing a single crossing, hence the unknotting number is one for both of them. Note that a knot has unknotting number zero if and only if it is the unknot.

There is a natural way to add together two knots $K_1$ and $K_2$, called the \emph{connected sum}, denoted $K_1 \# K_2$. Here is the connected sum of the trefoil and the Conway knot:
\begin{center}
\includegraphics[scale=0.8]{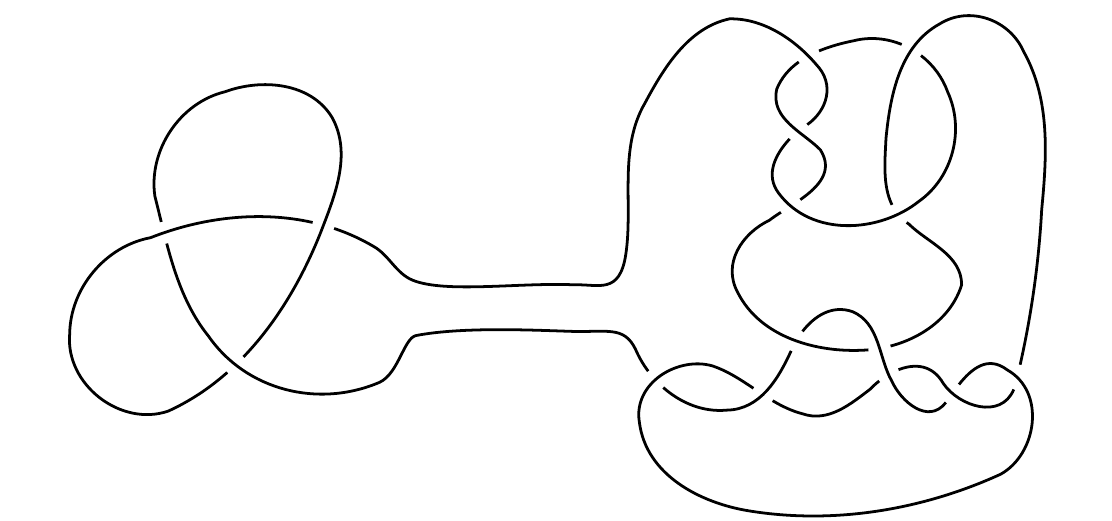} 
\end{center}

What is the unknotting number of $K_1 \# K_2$? A natural guess is that $u(K_1 \# K_2) = u(K_1) + u(K_2)$. One can readily check that $u(K_1 \# K_2) \leq u(K_1) + u(K_2)$. However, whether or not the reverse inequality holds remains an open question!

Here is another measure of complexity. Every knot in the 3-sphere bounds a compact, oriented, connected surface. Such surface is called a \emph{Seifert surface} for the knot. Recall that compact, oriented surfaces with connected boundary are characterized up to homeomorphism by their genus. The surfaces below are all have genus one:
\begin{center}
\includegraphics[scale=0.8]{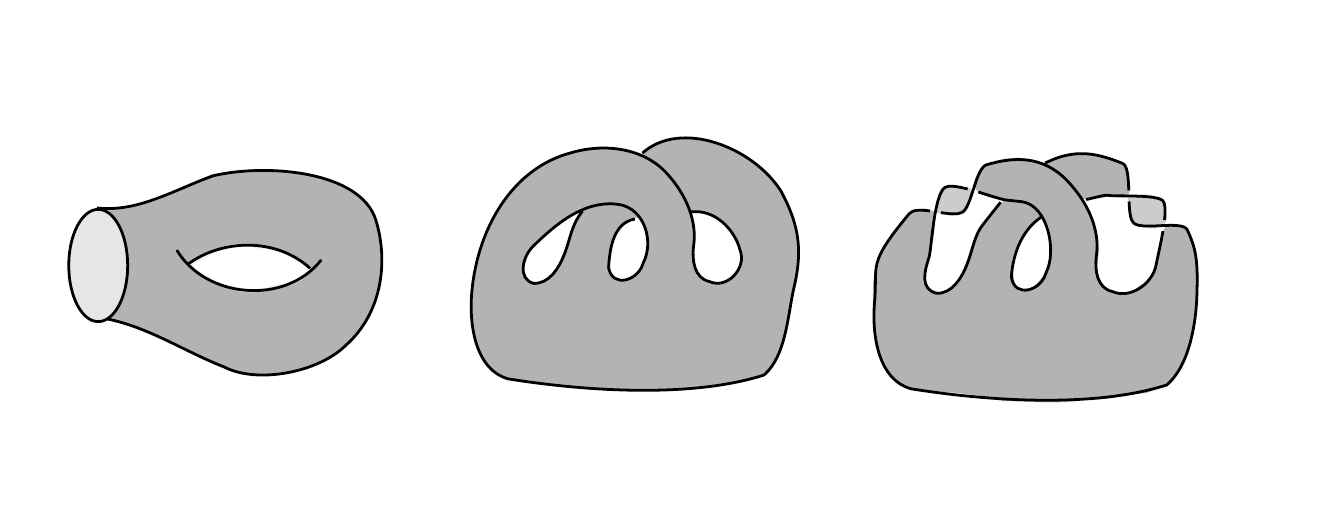} 
\end{center}
The boundary of each of the first two surfaces is the unknot. The boundary of the last surface is the trefoil.

The \emph{genus} of a knot $K$ is the minimal genus of a Seifert surface for $K$. The unknot is the only knot that bounds a disk. In other words, a knot had genus zero if and only if is the unknot. In contrast to unknotting number, we know how genus behaves under connected sum; Schubert \cite{Schubert1949} showed that genus is additive under connected sum, that is, $g(K_1 \# K_2) = g(K_1) + g(K_2)$.


The Alexander polynomial gives a lower bound on the genus of a knot: 
\[ \frac{1}{2} \deg \Delta_K(t) \leq g(K). \]
Since the Kinoshita-Terasaka knot and the Conway knot both have Alexander polynomial one, this bound does prove any useful information about their genera; for that, we turn to a result of Gabai, using foliations:
\begin{example}[\cite{Gabai}]
The Conway knot has genus three. The Kinoshita-Terasaka knot has genus two.
\end{example}

The unknot is the only knot with unknotting number zero, and it's also the only knot with genus zero. What about a measure of complexity where there are nontrivial knots that are also simple? Enter the \emph{slice genus}.

Recall that $S^3$ is the boundary of the 4-ball, and that a knot $K$ in $S^3$ is \emph{slice} if it bounds a smooth disk in the 4-ball. Such a disk is a called a \emph{slice disk} for $K$. Not every knot $K$ bounds a smooth disk in the 4-ball, but every knot does bound a smooth compact, oriented, connected surface in the 4-ball. (One way to obtain such a surface is by pushing a Seifert surface for $K$ into the 4-ball.) The minimal genus of such surface is called the \emph{slice genus} of $K$. Slice knots are precisely those knots with slice genus zero. Of course the unknot is slice, but there are also infinitely many nontrivial knots which are slice. For example, the Kinoshita-Terasaka knot is slice. Unlike the ordinary genus of a knot, slice genus is not additive under connected sum. 

The Alexander polynomial can obstruct sliceness: if $K$ is slice, then $\Delta_K(t)$ is of the form $t^n f(t) f(t^{-1})$ for some polynomial $f$ and some natural number $n$.

\begin{example}
The trefoil is not slice, since its Alexander polynomial $t^2-t+1$ is irreducible.
\end{example}

Closely related to the notion of sliceness is the following equivalence relation: two knots $K_0$ and $K_1$ are \emph{concordant} if they cobound an annulus $A$ in $S^3 \times [0,1]$, where the boundary of $A$ is $K_0 \subset S^3 \times \{ 0\}$ and $K_1 \subset S^3 \times \{ 1\}$. One can check that a knot is slice if and only if it concordant to the unknot.


Note that we required our surfaces to be smoothly embedded. What would happen if we just asked for topologically embedded disks in $B^4$? It turns out that every knot bounds a topologically embedded disk in $B^4$. Recall that the \emph{cone} of a space $X$ is $\Cone(X) = (X \times [0,1])/(X \times \{0\})$. Since $\Cone (S^3, K) = (B^4, B^2)$, every knot $K$ in $S^3$ bounds a topological disk in $B^4$, but the disk is not smoothly embedded, because of the cone point. Rather than requiring smoothness, one can instead require that the disk be locally flat; a knot that bounds a locally flat disk is called \emph{topologically slice}. Freedman \cite{Freedmandisk} proved that any knot with Alexander polynomial one is topologically slice; in particular, the Conway knot is topologically slice. Work of Donaldson \cite{Donaldson} implies that there are topologically slice knots that are not slice. Many slice obstructions actually obstruct topological sliceness, which is part of the reason why showing the Conway knot is not slice is so difficult.


\section{An equivalent condition for sliceness}
There are many invariants that obstruct sliceness, such as the aforementioned factoring of the Alexander polynomial, integer-valued invariants $\tau$ and $\nu$ coming from knot Floer homology \cite{OS4ball, OSrational}, and Rasmussen's integer-valued invariant $s$ coming from Lee's perturbation of Khovanov homology \cite{Rasmussen4ball, Lee}. These invariants (and many more!) all vanish for the Conway knot. 
(In my PhD thesis, I defined a new slice obstruction. One of the first questions people asked me was what its value was on the Conway knot; sadly, the obstruction vanishes for the Conway knot.)

Recall that in Section \ref{sec:tellingapart}, starting from a knot $K$ in $S^3$, we built a 3-manifold, the knot complement. Piccirillo's strategy for showing that the Conway knot is not slice relies on building a 4-manifold, called the \emph{knot trace}, from a knot $K$ in $S^3$. We will denote the trace of $K$ by $X(K)$. The following folklore result (see \cite{FoxMilnor}) is a key ingredient in Piccirillo's proof:

\begin{namedthm*}{Trace Embedding Lemma}
A knot $K$ is slice if and only if its trace $X(K)$ smoothly embeds in $S^4$.
\end{namedthm*}

In contrast to the fact that knots are determined by their complements, knots are not determined by their traces. That is, there exist non-isotopic knots $K$ and $K'$ with the same (i.e., diffeomorphic) traces \cite{Akbulut2dim}. 
Allison Miller and Lisa Piccirillo \cite{MillerPiccirillo} proved something even stronger: they showed that there exist knots $K$ and $K'$ with the same trace such that $K$ and $K'$ are not even concordant. This disproved a conjecture of Abe \cite{Abe}. Miller and Piccirillo's result implies that it's possible to have knots $K$ and $K'$ with the same trace, but for, say, $s(K)$ to be zero while $s(K')$ is nonzero.

We are slowing uncovering Piccirillo's strategy for proving the Conway knot is not slice: find a knot $K'$ with the same trace as the Conway knot, and show that $K'$ is not slice. Then the Trace Embedding Lemma implies that the Conway knot is not slice either.

\section{Handles and traces}

Let $B^n$ denote the $n$-ball.
Recall the 3-ball with a handle attached from beginning of these notes. More specifically, the handle consists of $B^1 \times B^2$ attached to $S^2 = \partial B^3$ along $S^0 \times B^3 = \partial B^1 \times B^2$. This handle is called a $3$-dimensional $1$-handle. 

\begin{center}
\labellist
	\pinlabel {attaching region} at 133 93
	\pinlabel {core} at 96 140
\endlabellist
\includegraphics[scale=1]{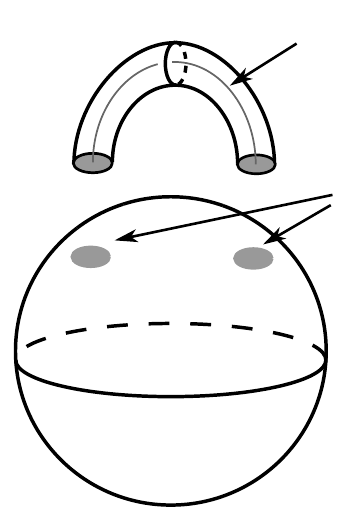} 
\end{center}

More generally, we consider an $n$-dimensional $k$-handle $B^k \times B^{n-k}$. Such a handle can be attached to an $n$-manifold $M$ with boundary by identifying a submanifold $S^{k-1} \times B^{n-k} \subset \partial M$ with $S^{k-1} \times B^{n-k} = \partial B^k \times B^{n-k}$. The submanifold $S^{k-1} \times B^{n-k} \subset \partial M$ is called the \emph{attaching region} of the handle. The \emph{core} of the handle is $B^k \times \{0\}$, where we think of $B^k$ as the unit ball in $\R^k$.



To build the knot trace, we will consider a $4$-dimensional $2$-handle $B^2 \times B^2$ attached to $S^3 = \partial B^4$. We need to specify the attaching region $S^1 \times B^2 \subset S^3$. This is just a tubular neighborhood of a knot. (The careful reader will note that we need to specify a parametrization of the neighborhood with $S^1 \times B^2$; this is called the \emph{framing} of the knot. For ease of exposition, we will largely suppress this key point from our discussion.) The \emph{trace} of a knot $K$ is the result of attaching a ($0$-framed) $2$-handle to $S^3 = \partial B^4$ along $K$. This is just a higher dimensional analog of the 1-handle attached to the 3-ball above.

\section{Knots with the same trace}
In order to understand Piccirillo's construction of a knot with the same trace as the Conway knot, it will be helpful to consider an analogy one dimension lower, in 3-dimensions, where we can more easily visualize things.

Consider the 3-ball with a 1-handle attached. Recall that a (3-dimensional) 2-handle is just a thickened disk $B^2 \times B^1$, which we attached along an annulus $S^1 \times B^1$. Suppose we attached a 2-handle along the grey annulus:

\begin{center}
\includegraphics[scale=1]{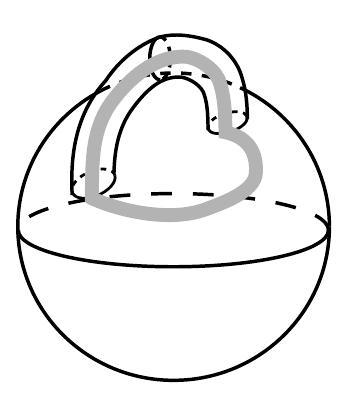} 
\end{center}
Observe that the resulting manifold $M_1$ is homeomorphic (in fact, diffeomorphic, after smoothing corners) to $B^3$! 

We could instead attach a 2-handle along the following grey thickened curve:
\begin{center}
\includegraphics[scale=1]{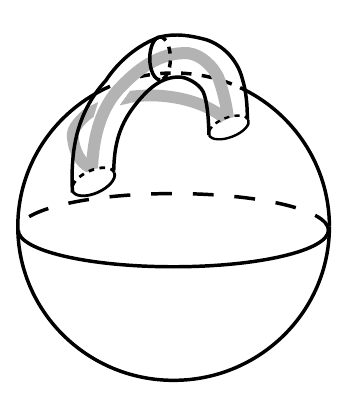} 
\end{center}
This would yield a manifold, $M_2$, which is again homeomorphic to $B^3$.

If we attached 2-handles to both of the grey curves, we obtain a manifold $M$ that is homeomorphic to $B^3$ with a 2-handle attached. Note that $M$ is built from a 3-ball, one 1-handle, and two 2-handles. We could view $M$ as $M_1 \cong B^3$ with a 2-handle attached 
or we could view $M$ as $M_2 \cong B^3$ with a 2-handle attached. 
Notice that the attaching regions for these 2-handles are just (thickened) embedded circles in $S^2 = \partial B^3$. Of course, embedded circles in $S^2$ are not especially interesting. But what happens when we bump things up a dimension?

Now consider the trace $X(C)$ of the Conway knot $C$. Piccirillo found a clever way to build $X(C)$ as a 4-ball, a 1-handle, and two 2-handles. (All of the handles here are 4-dimensional.) If you take the 4-ball, the 1-handle, and the first 2-handle, you get a 4-ball, and the second 2-handle is attached along the Conway knot $C$ in $S^3$ (the boundary of the 4-ball). On the other hand, if you take the 4-ball, the 1-handle, and the second 2-handle, you still get a 4-ball, and the remaining 2-handle is attached along some different knot $K'$. This means that $C$ and $K'$ have the same trace! Here is Piccirillo's knot $K'$ that has the same trace as the Conway knot:
\begin{center}
\includegraphics[scale=1]{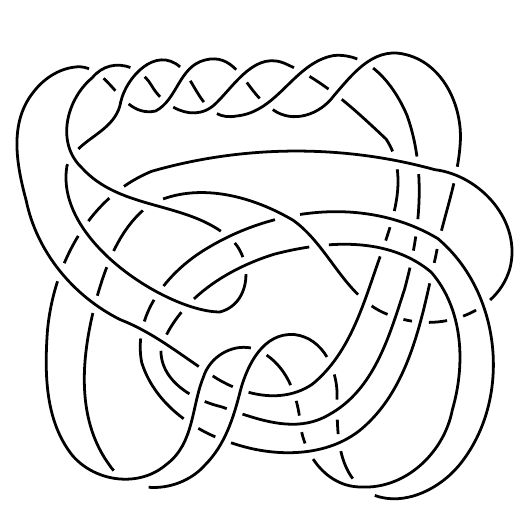} 
\end{center}

\section{Proof of the Trace Embedding Lemma}

Now that we have seen handles and traces, we will sketch the proof of the Trace Embedding Lemma. 

Suppose that $K$ is slice. This means that $K$ bounds a smooth disk in the 4-ball. Recall that $S^4$ is the union of two 4-balls, say $B^4_1$ and $B^4_2$. Think of $K$ as sitting in the common $S^3$ boundary of these two 4-balls. Since $K$ is slice, it bounds a slice disk $D$ in say $B^4_2$. Recall that a 4-dimensional 2-handle is just $D^2 \times D^2$. Then $B^4_1$ together with a closed neighborhood of $D$ is the trace of $K$, smoothly embedded in $S^4$. A schematic of $S^4$ as the union of two 4-balls is shown below:
\begin{center}
\labellist
	\pinlabel {$B_1^4$} at 76 40
	\pinlabel {$B_2^4$} at 76 120
	\pinlabel {$K$} at 72.5 82
	\pinlabel $S^3$ at 25 70
\endlabellist
\includegraphics[scale=1]{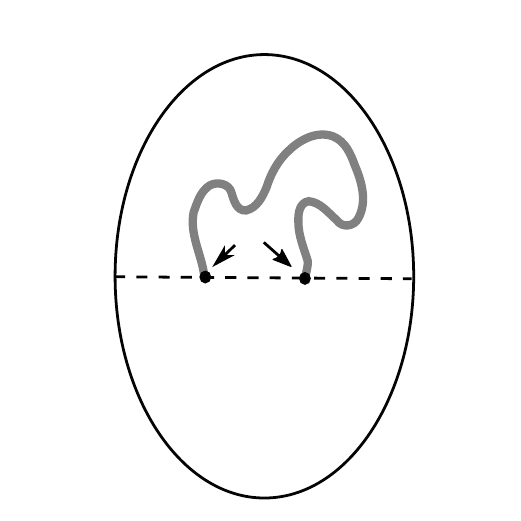}
\end{center}
The slice disk is represented by the thick grey curve. The trace of $K$ consists of $B_1^4$ together with a neighborhood of the slice disk for $K$.


Now suppose that $X(K)$ embeds in $S^4$. Consider the piecewise linear embedded $S^2$ in $X(K)$ consisting of the core of the $2$-handle together with the cone of $K$. Smoothly embed $X(K)$ in $S^4$; composition gives a piecewise linear embedding of $S^2$ in $S^4$, which is smooth away from the cone point $p$. Now take a small neighborhood around $p$ in $S^4$. The complement of this neighborhood is a 4-ball $B$. Consider the piecewise linear embedding of $S^2$ intersected with $B$; we've cut out the cone point, so this gives a slice disk in $B$ for $K$ in $\partial B$. A schematic of the trace embedded in $S^4$ is shown below:

\begin{center}
\labellist
	\pinlabel {\small{$p$}} at 68 50
	\pinlabel {$S^3$} at 85 45
	\pinlabel {slice disk for $K$} at -10 50
	\pinlabel {$K$} at 75 80
\endlabellist
\includegraphics[scale=1]{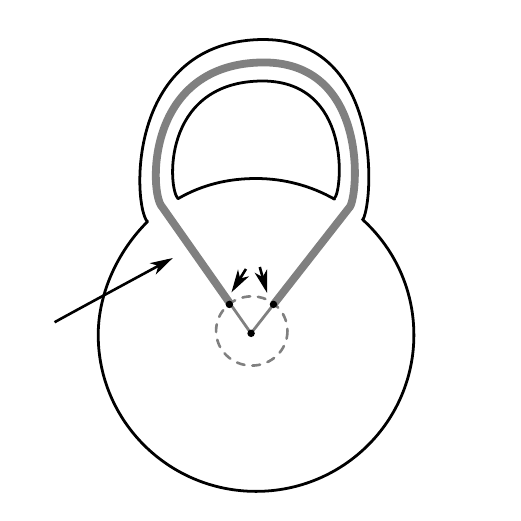}
\end{center}
The 4-ball $B$ is everything outside of the $S^3$ dotted circle, and the thick grey curve shows the slice disk for $K$.

\section{Showing that $K'$ is not slice}
The goal is now to find a way to show that $K'$, the knot that shares a trace with the Conway knot, is not slice. It turns out that some slice obstructions, such as the invariant $\nu$ coming from knot Floer homology, are actually trace invariants: if two knots $K_1$ and $K_2$ have the same trace, then $\nu(K_1) = \nu(K_2)$ \cite{HMP}. 

Luckily, the same is not true for Rasmussen's $s$-invariant. Using a computer program and some simple algebraic observations, Piccirillo shows that $s(K')=2$, implying that $K'$ is not slice. Since $K'$ and the Conway knot have the same trace, the Trace Embedding Lemma implies that the Conway knot is not slice. 

\section{What's next?}
Now that we know exactly which knots with fewer than 13 crossings are slice, what's next? Of course, one could try to determine exactly which knots with fewer than 14 or 15 crossings are slice. But why not try to apply some of our tools to other open problems? 

The smooth 4-dimensional Poincar\'e conjecture posits that a smooth 4-manifold 
that is homeomorphic to $S^4$ is actually diffeomorphic to $S^4$. To disprove the conjecture, one wants to find an \emph{exotic} $S^4$, that is, a smooth 4-manifold that is homeomorphic but not diffeomorphic to $S^4$.
One possible approach (outlined in 
\cite{FGMW}) to disprove the smooth 4-dimensional Poincar\'e conjecture relies on Rasmussen's $s$-invariant, as follows.

There are many constructions of potentially exotic 4-spheres $\Sigma$ (see, for example \cite{CS}; note that certain infinite subfamilies of these are known to be standard by \cite{Akbulut, GompfCS, MeierZupan}). 
By removing a neighborhood of a point in $\Sigma$, one can instead study potentially exotic 4-balls $\beta$. The difficult part is now determining whether or not $\beta$ is exotic, or if it is in fact just the standard $B^4$.

While slice obstructions like $\nu$ actually obstruct a knot from being slice in an exotic $4$-ball, it remains possible that the $s$-invariant only obstructs a knot from being slice in the standard 4-ball. The game is then to try to find a knot $K$ that is slice in a potentially exotic 4-ball $\beta$. If $s(K)$ is non-zero, then $K$ is not slice in the standard 4-ball, thereby implying that $\beta$ must be exotic.

Both of the key steps in this approach (constructing the potentially exotic 4-ball and computing $s$) seem difficult. But maybe there is some other way to get handle on the problem in order to trace a solution. I look forward to reading an article about such a result!

\section*{Acknowledgements}
I would like to thank JungHwan Park and Lisa Piccirillo for helpful comments on an earlier draft.

\bibliographystyle{amsalpha}
\def\MR#1{}
\bibliography{bib}

\end{document}